\documentclass[12pt]{amsart}
\usepackage{amsmath}
\usepackage{amsthm}
\usepackage{amssymb}
\usepackage{color}
\usepackage{indentfirst}
\usepackage{graphicx}
\usepackage{bbding}
\usepackage{hyperref}

\begin{document}

\newtheorem{theorem}{Theorem}
\newtheorem{lemma}{Lemma}
\newtheorem{corollary}{Corollary}
\newtheorem{proposition}{Proposition}
\newtheorem*{Ithm}{Ibragimov's Theorem}

\theoremstyle{definition}
\newtheorem{definition}{Definition}
\newtheorem{example}{Example}
\newtheorem{remark}{Remark}

\newenvironment{prof}[1][Proof]{\noindent\textit{#1}\quad }

\newcommand{\ind}{\mathop{\mathrm{index}}}
\parindent2em
\title{On the quasisymmetric minimality of homogeneous perfect sets}
\author{Yingqing Xiao }
\address{Yingqing Xiao, College of Mathematics and Econometrics, Hunan University, Changsha, 410082, China}
\email{ouxyq@hnu.edu.cn}

\author{Zhanqi Zhang}
\address{Zhanqi Zhang, College of Mathematics and Econometrics, Hunan University, Changsha, 410082, China}
\email{rateriver@sina.com}
\maketitle

\textbf{Abstract:} Z. Wen and J. Wu introduced the notion of homogeneous perfect sets as a generalization of Cantor type sets and determined their exact Hausdorff dimension based on the length of their basic intervals and the gaps between them. In this paper, we considered the minimality of the homogeneous perfect sets with Hausdorff dimension $1$ and proved they are $1$-dimensional quasisymmetrically minimal under some conditions.


\smallskip
\textbf{2000 mathematics classification:} Primary 30C62; Secondary 28A78.

\section{Introduction}
Let $(X,d_{X})$, $(Y,d_{Y})$ be two metric spaces. A topological homeomorphism $f:X\rightarrow Y$ is called quasisymmetric if there is a homeomorphism $\eta:[0,\infty)\rightarrow[0,\infty)$ such that
\begin{equation}
\ \frac{d_{Y}(f(x),f(a))}{d_{Y}(f(x),f(b))}\leq\eta(\frac{d_{X}(x,a)}{d_{X}(x,b)})\nonumber
\end{equation}
for all triples $a,b,x$ of distinct points in $X$. In particular, we also say that $f$ is a $n$-dimensional quasisymmetric mapping when $X=Y=\mathbb{R}^{n}$.

We call a set $E\subseteq \mathbb{R}^{n}$ $n$-dimensional quasisymmetrically minimal if
\begin{equation}
\ \dim_{H}f(E)\geq \dim_{H}E\nonumber
\end{equation}
for any $n$-dimensional quasisymmetric mapping $f$, where $\dim_{H}$ denotes the Hausdorff dimension. For convenience, we call a set minimal if it is $n$-dimensional quasisymmetrically minimal. When $n\geq2$, all sets in $\mathbb{R}^{n}$ with Hausdorff dimension $n$ are minimal \cite{G}. But according to \cite{T}, a set in $\mathbb{R}$ with positive Lebesgue measure may not be  minimal. So if we have a set in $\mathbb{R}$ with Hausdorff dimension 1, it is interesting to ask whether it is minimal or not.
 \\
\indent Obviously, sets with nonempty interior are minimal. In \cite{SW}, the authors introduced quasisymmetrically thick sets which are minimal. Recently, Wang and Wen proved that uniform Cantor sets with Hausdorff dimension $1$ are minimal in \cite{WWen,WW}. Besides, a large class of Moran sets of Hausdorff dimension 1 are minimal \cite{D}.

In this paper, we prove that homogeneous perfect sets with Hausdorff dimension $1$ are minimal under some conditions. In \cite{W}, the notion of homogeneous perfect sets as a generalization of uniform Cantor sets was introduced. Uniform Cantor sets are special cases of homogeneous perfect sets, so our paper generalized the results of \cite{WWen,WW}. On the other hand, Moran sets are more general than homogeneous perfect sets. But in \cite{D}, the authors proved that Moran sets with Hausdorff dimension $1$ are minimal under some conditions including that the numbers of basic intervals of previous one are bounded. In order to remove this condition, we add some restrictions on the size of gaps and the lengths of basic intervals, which turn out to be the homogeneous perfect sets case.
\\
\indent We organize the paper as follows: In Section 2, we introduce homogeneous perfect sets and show our main results. After that we do some reconstructions to the sets which is necessary for the proof; When studying quasisymmetric minimality of Cantor type sets, we have found some basic routine to check its minimality. In Section 3, we introduce the routine; Then all we have to do is to check some conditions of homogeneous perfect sets, we do these in Section 4; In the last section, we introduce an example in which the homogeneous perfect set does not satisfy some condition we assumed in Theorem \ref{theorem-1}, but it is still minimal.

\section{Homogeneous perfect sets and reconstructions}
\textbf{Hausdorff dimension}.
Let $K\subseteq \mathbb{R}^n$. For any $s\geq 0$, the $s$-dimensional Hausdorff measure of $K$ is given in the usual way by
$$
\mathbf{H}^s(K)=\liminf_{\delta\rightarrow 0}\{\sum_i|U_i|^s: K\subset \bigcup_iU_i, 0<|U_i|<\delta\}.
$$
This leads to the definition of the
Hausdorff dimension of $K$:
$$
\dim_HK=\inf\{s: \mathbf{H}^s(K)<\infty\}=\sup\{s: \mathbf{H}^s(K)>0\}.
$$
For more details on Hausdorff dimension we refer to \cite{F}.

Now, we recall the notion of homogeneous perfect sets.
\subsection{Homogeneous perfect sets}
Let $J_{\emptyset}=[0,1]$. Let $\{n_{k}\}_{k\geq1}$ be a sequence of positive integers and $\{c_{k}\}_{k\geq1}$ a sequence of positive real numbers such that for any $k\geq1$, $n_{k}\geq2$ and $0<c_{k}<1$. For any $k\geq1$, let $\Omega_{k}=\{(i_{1},i_{2},\cdots,i_{k}):1\leq i_{j}\leq n_{j},1\leq j\leq k\}$, $\Omega=\bigcup_{k\geq 0}\Omega_{k}$, where $\Omega_{0}=\{\emptyset\}$. If $\sigma=(\sigma_{1},\sigma_{2},\cdots,\sigma_{k})\in\Omega_{k}$, $1\leq j\leq n_{k+1}$, then $\sigma\ast j:=(\sigma_{1},\sigma_{2},\cdots,\sigma_{k},j)\in\Omega_{k+1}$. Let $\mathcal{J}=\{J_{\sigma}:\sigma\in\Omega\}$ be a collection of closed subintervals of $J_{\emptyset}$. We say that the collection $\mathcal{J}$ fulfill the homogeneous perfect structure provided:
\\
\indent 1. For any $k\geq0$ and $\sigma\in\Omega_{k}$, $J_{\sigma\ast 1},J_{\sigma\ast 2},\cdots,J_{\sigma\ast n_{k+1}}$ are subintervals of $J_{\sigma}$. Furthermore, $\max\{x:x\in J_{\sigma\ast i}\}\leq \min\{x:x\in J_{\sigma\ast (i+1)}\}$ for $1\leq i\leq n_{k+1}-1$, that is the interval $J_{\sigma\ast i}$ is located at the left of $J_{\sigma\ast (i+1)}$ and the interiors of the intervals $J_{\sigma\ast i}$ and $J_{\sigma\ast (i+1)}$ are disjoint;
\\
\indent 2. For any $k\geq 1$, $\sigma\in\Omega_{k-1}$ and $1\leq i\leq n_{k}$, we have
\begin{equation}
\ \frac{|J_{\sigma\ast i}|}{|J_{\sigma}|}=c_{k}.\nonumber
\end{equation}
Here and in the sequel $|\cdot|$ stands for the 1-dimensional Lebesgue measure;
\\
\indent 3. There exists a sequence of nonnegative real numbers $\{\eta_{k,j}:k\geq 1,0\leq j\leq n_{k}\}$ such that for any $k\geq0$ and $\sigma\in\Omega_{k}$, we have $\min(J_{\sigma\ast 1})-\min(J_{\sigma})=\eta_{k+1,0}$, $\max(J_{\sigma})-\max(J_{\sigma\ast n_{k+1} })=\eta_{k+1,n_{k+1}}$, and for any $1\leq i\leq n_{k+1}-1$, we have $\min(J_{\sigma\ast (i+1)})-\max(J_{\sigma\ast i})=\eta_{k+1,i}$.
\\
\indent Suppose that the collection of intervals $\mathcal{J}=\{J_{\sigma}:\sigma\in\Omega\}$ satisfies the homogeneous perfect structure, let $E_{k}=\bigcup_{\sigma\in\Omega_{k}}J_{\sigma}$ for any $k\geq 0$, then the set
\begin{equation}
\ E:=E(J_{\emptyset},\{n_{k}\},\{c_{k}\},\{\eta_{k,j}\})=\bigcap_{k\geq 0}E_{k} \nonumber
\end{equation}
is called a homogeneous perfect set. The intervals $J_{\sigma}$, $\sigma\in\Omega_{k}$ are called basic intervals of order $k$.
The Hausdorff dimension of homogeneous perfect set $E$ which depends on $\{n_k\}$, $\{c_k\}$ and $\{\eta_{k,j}\}$ have been obtained in \cite{W} as follows.
\begin{theorem}[\cite{W}, Theorem 1.2]\label{Hdim}
Let $E = E(J_{\emptyset},\{n_k\},\{c_k\},\{\eta_{k,j}\})$ be a homogeneous perfect set. Suppose
there exist positive constants $\tilde{c}_{1}$, $\tilde{c}_{2}$, $\tilde{c}_{3}$ and $\tilde{c}_{4}$ such that for any $k\geq1$, at least one of the following four conditions is satisfied:
\\
\indent $(1)$ $\max_{1\leq l\leq n_{k}-1}\eta_{k,l}\leq \tilde{c}_{1}\min_{1\leq l\leq n_{k}-1}\eta_{k,l}$;
\\
\indent $(2)$ $\max_{1\leq l\leq n_{k}-1}\eta_{k,l}\leq \tilde{c}_{2}\cdot c_{1}c_{2}\cdots c_{k}$;
\\
\indent $(3)$ $n_{k}\cdot\min_{1\leq l\leq n_{k}-1}\eta_{k,l}\geq \tilde{c}_{3}\cdot c_{1}c_{2}\cdots c_{k-1}$;
\\
\indent $(4)$ $n_{k}\leq \tilde{c}_{4}$.
\\Then
\begin{equation}
\dim_HE=\liminf_{k\rightarrow \infty}\frac{\log(n_{1}n_{2}\cdot\cdot\cdot n_{k})}{-\log(\sum_{l=1}^{n_{k+1}-1}\eta_{k+1,l}+n_{k+1}c_1c_2\cdot\cdot\cdot c_{k+1})}.\nonumber
\end{equation}
\end{theorem}

In \cite{X}, Xiao proved homogeneous perfect sets with Hausdorff dimension $1$ are minimal under the condition $(4)$ of Theorem \ref{Hdim}. In the following, we study the minimality of the homogeneous perfect sets under other three conditions. After reconstruction of homogeneous perfect sets, we will prove our main result:
\begin{theorem}\label{theorem-1}
Suppose $E=E(J_{\emptyset},\{n_{k}\},\{c_{k}\},\{\eta_{k,j}\})$ is a homogeneous perfect set with $\dim_{H}E=1$, then $E$ is minimal if it satisfies the condition $(A)$ or $(B)$:

$(A)$ there exists $C\geq 1$ such that $\max_{1\leq l\leq n_{k}-1}\eta_{k,l}\leq C \min_{1\leq l\leq n_{k}-1}\eta_{k,l}$ for all $k\geq 1$;

$(B)$ there exists $D>0$ such that $\max_{1\leq l\leq n_{k}-1}\eta_{k,l}\leq D c_{1}c_{2}\cdots c_{k}$ for all $k\geq 1$.
\end{theorem}

Unfortunately, we cannot prove that homogeneous perfect sets with Hausdorff dimension $1$ are minimal under the condition $(3)$ of Theorem \ref{Hdim}.
\\
\indent Recently, Yang also considered the minimality of the homogeneous perfect sets and obtained the following result in \cite{Yang}.
\begin{theorem}[\cite{Yang}, Theorem 1]\label{Y-Hdim}
Let $E = E(J_{\emptyset},\{n_k\},\{c_k\},\{\eta_{k,j}\})$ be a homogeneous perfect set satisfying that there are $L\geq 1$, $c>0$ such that $\frac{1}{L}\leq\frac{\eta_{k,i}}{\eta_{k,j}}\leq L$ and $\eta_{k,0}+\eta_{k,n_k}\leq c\eta_{k,i}$ for all $k\geq 1, 1\leq i,j\leq n_k -1$. If $\mathrm{dim}_HE=1$, then $\mathrm{dim}_H f(E)=1$ for any $1$-dimensional quasisymmetric mapping $f$.
\end{theorem}
Obviously, the condition $(A)$ of our Theorem \ref{theorem-1} is equivalent to the condition: there is $L\geq 1$ such that $\frac{1}{L}\leq\frac{\eta_{k,i}}{\eta_{k,j}}\leq L$ for all $k\geq 1$, $1\leq i,j\leq n_k -1$.

\subsection{Reconstruction of homogeneous perfect sets}
Next we do some reconstructions. We did step 1 in order to weaken the influence of leftmost gaps and rightmost gaps in every basic intervals and to get better use of Theorem \ref{Hdim}. We did step 2 in order to deal with the condition that $\sup_k\{n_{k}\}=\infty$. If $\sup_k\{n_{k}\}<\infty$, we can still take this step. The idea of step 2 is based on \cite{WW}.
\\
\indent \textbf{Step 1}: For any $k\geq 0$ and $\sigma\in\Omega_{k}$, let $J_{\sigma}^{\ast}$ be a closed subinterval of $J_{\sigma}$ such that
\\
\indent (1) the distance between the left endpoint of $J_{\sigma}^{\ast}$ and the left endpoint of $J_{\sigma}$ is $\eta_{k+1,0}$;
\\
\indent (2) the distance between the right endpoint of $J_{\sigma}^{\ast}$ and the right endpoint of $J_{\sigma}$ is $\eta_{k+1,n_{k+1}}$. So we have $|J_{\sigma}^{\ast}|=\Sigma_{l=1}^{n_{k+1}-1}\eta_{k+1,l}+n_{k+1}c_{1}c_{2}\cdots c_{k+1}$. We denote $\delta_{k}:=|J_{\sigma}^{\ast}|$ for $k\geq1$ and $\delta_{0}:=|J_{\emptyset}^{\ast}|$.
\\
\indent Let $E_{k}^{\ast}=\bigcup_{\sigma\in\Omega_{k}}J_{\sigma}^{\ast}$ for all $k\geq 0$, we have
\begin{equation}
\ E=\bigcap_{k\geq 0}E_{k}=\bigcap_{k\geq 0}E_{k}^{\ast}.\nonumber
\end{equation}
In fact, $E=E(J_{\emptyset}^{\ast},\{n_{k}\},\{c_{k}^{\ast}\},\{\eta_{k,j}^{\ast}\})$ is a homogenous perfect set with the following parameters:

$(a)$ $J_{\emptyset}^{\ast}=[0,1]-[0,\eta_{1,0})-(\eta_{1,n_{1}},1]$;

$(b)$ $c_{k}^{\ast}=\frac{\delta_{k}}{\delta_{k-1}}$;

$(c)$ $\eta_{k,l}^{\ast}=\eta_{k+1,0}+\eta_{k+1,n_{k+1}}+\eta_{k,l}$ for $1\leq l\leq n_{k}-1$ and $\eta_{k,0}^{\ast}=\eta_{k+1,0}$, $\eta_{k,n_{k}}^{\ast}=\eta_{k+1,n_{k+1}}$.

From above, in order to get homogeneous perfect set $E$, we find a sequence of sets $\{E_{k}^{\ast}\}_{k\geq0}$. For each $E_{k}^{\ast}$, it has basic intervals $J_{\sigma}^{\ast}$ , $\sigma\in\Omega_{k}$ with length $\delta_{k}$. For each $J_{\sigma}^{\ast}$, $J_{\sigma}^{\ast}\bigcap E_{k+1}^{\ast}$ has $n_{k+1}$ basic intervals $J_{\sigma\ast l}^{\ast}$ , $1\leq l\leq n_{k+1}$ with equal length $\delta_{k+1}$. And the gaps between them in $J_{\sigma}^{\ast}$ have lengths $\eta_{k+1,0}^{\ast},\eta_{k+1,1}^{\ast},\cdots,\eta_{k+1,n_{k+1}}^{\ast}$ from left to right.
\\
\indent \textbf{Step 2}: In order to get the same homogeneous perfect set $E$, we try to find another sequence of sets $\{F_{m}\}_{m\geq0}$ such that $\{E_{k}^{\ast}\}_{k\geq0}$ is a subsequence of $\{F_{m}\}_{m\geq0}$ and $E=\bigcap_{m\geq0}F_{m}$.
\\
\indent For every $k\geq 1$, let $i_{k}$ be a positive integer such that $2^{i_{k}}\leq n_{k}<2^{i_{k}+1}$. Let $m_{0}=0$ and $m_{k}=i_{1}+i_{2}+\cdots+i_{k}$ for $k\geq 1$. Set $F_{m_{k}}=E_{k}^{*}$ for all $k\geq 0$. Given $k\geq 1$, we are going to construct $F_{m}$ for $m_{k-1}<m<m_{k}$. If $2\leq n_{k}\leq 3$, $i_{k}=1$ and nothing need to do. If $n_{k}\geq 4$, we are to construct $F_{m_{k-1}+1}$ first.
\\
\indent Construct $F_{m_{k-1}+1}$: Write $n_{k}=a_{0}+a_{1}\cdot2+a_{2}\cdot2^{2}+\cdots+a_{i_{k}-1}\cdot2^{i_{k}-1}+2^{i_{k}}$, where $a_{j}=0$ or $1$ for $0\leq j\leq i_{k}-1$. For any $\sigma\in\Omega_{k-1}$, $J_{\sigma}^{\ast}$ have $n_{k}$ children basic intervals $J_{\sigma\ast1}^{\ast},\cdots,J_{\sigma\ast n_{k}}^{\ast}$ in $E_{k}^{\ast}$. Let $l_{1}=a_{1}+a_{2}\cdot2+\cdots+a_{i_{k}-1}\cdot2^{i_{k}-2}+2^{i_{k}-1}$ and let $J_{1}^{\sigma,1}=[J_{\sigma\ast1}^{\ast},J_{\sigma\ast2}^{\ast},\cdots,J_{\sigma\ast l_{1}}^{\ast}]$, $J_{2}^{\sigma,1}=[J_{\sigma\ast(l_{1}+1)}^{\ast},\cdots,J_{\sigma\ast n_{k}}^{\ast}]$. Here and in the sequel, for any adjacent intervals $I_{1},I_{2},\cdots,I_{t}$, $[I_{1},I_{2},\cdots,I_{t}]$ denotes the smallest closed interval containing them. In the end, we put $F_{m_{k-1}+1}=\bigcup_{\sigma\in\Omega_{k-1}}(J_{1}^{\sigma,1}\bigcup J_{2}^{\sigma,1})$. If $i_{k}=2$, we have done all. Otherwise, if $i_{k}\geq3$, we are going to construct $F_{m_{k-1}+2}$.
\\
\indent Construct $F_{m_{k-1}+2}$: Let $l_{2}=a_{2}+a_{3}\cdot2+\cdots+a_{i_{k}-1}\cdot2^{i_{k}-3}+2^{i_{k}-2}$ and let $a=\min\{a_{0},a_{1}\}$. For each $\sigma\in\Omega_{k-1}$, we get $2$ children basic intervals $J_{1}^{\sigma,1}, J_{2}^{\sigma,1}$ in $F_{m_{k-1}+1}$. For each one, it may be seen as the union of some basic intervals of $E_{k}^{\ast}$ and the gaps between them. Remove the midmost gap and we get two closed subintervals. It is worth noting that we may have some gaps to be empty, i.e., $\eta_{k,l}^{\ast}=0$ for some $l$, but we still consider it as a gap. When talking about the position of the gaps (to find the midmost gap), the empty gaps still count. If it has two midmost gaps, remove the left one. Through this, for each $\sigma\in\Omega_{k-1}$, we get $4$ closed subintervals of $J_{\sigma}^{\ast}$, that is
\\
\indent$J_{1}^{\sigma,2}=[J_{\sigma\ast1}^{\ast},\cdots,J_{\sigma\ast l_{2}}^{\ast}],J_{2}^{\sigma,2}=[J_{\sigma\ast(l_{2}+1)}^{\ast},\cdots,J_{\sigma\ast l_{1}}^{\ast}],
\\
\indent J_{3}^{\sigma,2}=[J_{\sigma\ast(l_{1}+1)}^{\ast},\cdots,J_{\sigma\ast (l_{1}+l_{2}+a)}^{\ast}],J_{4}^{\sigma,2}=[J_{\sigma\ast(l_{1}+l_{2}+a+1)}^{\ast},\cdots,J_{\sigma\ast n_{k}}^{\ast}]$.\\In the end, we put $F_{m_{k-1}+2}=\bigcup_{\sigma\in\Omega_{k-1}}(J_{1}^{\sigma,2}\bigcup J_{2}^{\sigma,2}\bigcup J_{3}^{\sigma,2}\bigcup J_{4}^{\sigma,2})$. If $i_{k}=3$, we have done all.

Construct $F_{m_{k-1}+j}$: We construct $F_{m_{k-1}+j}$ inductively when $1< j\leq i_{k}-1$. For any basic intervals of $F_{m_{k-1}+j-1}$, it may be seen as the union of some basic intervals of $E_{k}^{\ast}$ and the gaps between them. Remove the midmost gap and we can get two closed subintervals to be the basic intervals of $F_{m_{k-1}+j}$. Dealing with all basic intervals of $F_{m_{k-1}+j-1}$, we get children basic intervals of $F_{m_{k-1}+j}$. Uniting them all, we obtain $F_{m_{k-1}+j}$.
\\
\indent Construct $F_{m_{k-1}+i_{k}-1}$: We may continue process to get $F_{m_{k-1}+i_{k}-1}$. For each $\sigma\in\Omega_{k-1}$, we have $2^{i_{k}-1}$ basic intervals of $F_{m_{k-1}+i_{k}-1}$ and denote them by $J_{1}^{\sigma,i_{k}-1},\cdots,J_{2^{i_{k}-1}}^{\sigma,i_{k}-1}$. That is we have $F_{m_{k-1}+i_{k}-1}=\bigcup_{\sigma\in\Omega_{k-1}}(J_{1}^{\sigma,i_{k}-1}\bigcup\cdots\bigcup J_{2^{i_{k}-1}}^{\sigma,i_{k}-1})$.
\\
\indent Through above, we have found the sequence of sets $\{F_{m}\}_{m\geq0}$, it has properties:
\\
\indent (1) $F_{m+1}\subseteq F_{m}$ for $m\geq0$ and $E=\bigcap_{m\geq0}F_{m}$;
\\
\indent (2) For each basic interval of $F_{m}$, say $I_{m}$, it may have at most $4$ children basic intervals in $I_{m}\bigcap F_{m+1}$. Actually, if $n_{k}=a_{0}+a_{1}\cdot2+a_{2}\cdot2^{2}+\cdots+a_{i_{k}-1}\cdot2^{i_{k}-1}+2^{i_{k}}$ and $a_{i_{k}-1}=1$, we may find a basic interval of $F_{m_{k-1}+i_{k}-1}$ which have $4$ children basic intervals in $F_{m_{k}}$.
\\
\indent We have completed the reconstruction.

\section{Basic routine}
Before exhibiting the basic routine, we introduce some notations. Let $\mathcal{F}_{m}$ denote the family of all basic intervals of $F_{m}$. Let $\mathcal{G}_{m}$ denote all component intervals of $I_{m}-(I_{m}\bigcap F_{m+1})$, where $I_{m}\in\mathcal{F}_{m}$.
 That is $\mathcal{G}_{m}=\{$component intervals of $I_{m}-(I_{m}\bigcap F_{m+1}):I_{m}\in\mathcal{F}_{m}\}$.
Set $\mathcal{G}(I)=\{L: L\subset I, L\in \mathcal{G}_m\}$ for every $I\in \mathcal{F}_m$. Since for each basic interval of $F_{m}$, say $I$, it has at most $4$ children basic intervals in $I\bigcap F_{m+1}$, we have $\#(\mathcal{G}(I))\leq 5$. For $J\in\mathcal{F}_{m}$, if $m\geq1$, let $Fa(J)\in\mathcal{F}_{m-1}$ be the interval such that $J\subseteq Fa(J)$; If $m\geq0$, assuming that $J_1, J_2,\cdots, J_{N(J)}$ are the members of $\mathcal{F}_{m+1}$ which belong to $J$, we denote $C(J,1)=|J_1|+|J_2|+\cdots+|J_{N(J)}|$. Let
$$
\beta_{m}=\max\{\frac{|L|}{|I|}:I\in\mathcal{F}_{m} ,L\in\mathcal{G}_{m}, L\subseteq I\},
$$
$$
\Gamma_m=\min\{\frac{C(J,1)}{|J|}:J\in\mathcal{F}_{m}\}, \gamma_m=\max\{\frac{|J|}{|Fa(J)|}:J\in\mathcal{F}_{m}\}.
$$
Obviously, $0<\Gamma_m,\gamma_m\leq 1$ for all $m$. Since $\#(\mathcal{G}(I))\leq 5$, we have
\begin{lemma}\label{key-lemma} $\Gamma_{m}\geq 1-5\beta_m$ for all $m$.
\end{lemma}
\begin{prof}
For every $J
\in \mathcal{F}_m$, we have $\beta_m\geq \frac{|L|}{|J|}$, where $L\in \mathcal{G}(J)$. Thus, we have
$$
\sum_{L\in G(J)}\frac{|L|}{|J|}\leq \sum_{L\in G(J)}\beta_m\leq5\beta_m,
$$
which yields to
$$
\frac{C(J,1)}{|J|}=\frac{|J|-\sum_{L\in G(J)}|L|}{|J|}\geq 1-5\beta_m.
$$
Thus
$$
\Gamma_m=\min\{\frac{C(J,1)}{|J|}:J\in\mathcal{F}_{m}\}\geq 1-5\beta_m.
$$
\qed
\end{prof}
 After above introduction, we state the basic routine.
\begin{theorem}\label{theorem-2}
Let $E=\bigcap_{m\geq0}F_{m}$ be a homogeneous perfect set, then $E$ is minimal if it satisfies the conditions $(\tilde{a})$, $(\tilde{b})$, $(\tilde{c})$:

$(\tilde{a})$ $\lim\limits_{m\rightarrow\infty}\frac{1}{m}\sum_{j=0}^{m-1}\beta_{j}=0$;

$(\tilde{b})$ $\lim\limits_{m\rightarrow\infty}\frac{1}{m}\sum_{j=0}^{m-1}\log\Gamma_j=0$;

$(\tilde{c})$ there exists $\alpha\in(0,1)$ such that $\liminf\limits_{m\rightarrow\infty}\frac{\#\{1\leq i\leq m:\,\gamma_i<\alpha\}}{m}>0$.
\end{theorem}

Actually if we have a homogeneous perfect set $E$ with $\dim_{H}E=1$, in order to prove its minimality, all we have to do is to show that $\dim_{H}f(E)=1$ for any $1$-dimensional quasisymmetric mapping $f$. In this paper, the properties we need for quasisymmetric mappings are obtained in the next lemma.
\begin{lemma}[\cite{WWen}, Lemma 3.1]\label{Lemma-1}
Let $f$ be a $1$-dimensional quasisymmetric mapping. Then
\begin{equation}
\ \lambda\frac{|J|^{q}}{|I|^{q}}\leq\frac{|f(J)|}{|f(I)|}\leq4\frac{|J|^{p}}{|I|^{p}}\nonumber
\end{equation}
for all intervals $I,J$ with $J\subseteq I$, where $\lambda,p,q$ are three constants dependent on $f$ with $\lambda>0,0<p\leq 1\leq q$.
\end{lemma}

Fix any quasisymmetric mapping $f$ and let the meanings of $\lambda,p,q$ unchanged in the rest of this paper.
When we assume that $\dim_{H}E=1$, the information we get from this property is introduced in the next lemma.
\begin{lemma}\label{Lemma-2}
Let $E=E(J_{\emptyset},\{n_{k}\},\{c_{k}\},\{\eta_{k,j}\})$ be a homogeneous perfect set and $\dim_{H}E=1$. If the condition $(A)$ or $(B)$ in Theorem \ref{theorem-1} is satisfied, then we have
\begin{equation}
\ \lim\limits_{k\rightarrow\infty}\frac{\log_{2}(n_{1}n_{2}\cdots n_{k})}{-\log_{2}(\sum\limits_{l=1}^{n_{k+1}-1}\eta_{k+1,l}+n_{k+1}c_{1}c_{2}\cdots c_{k+1})}=1.\nonumber
\end{equation}
\end{lemma}
\begin{prof}
It is obvious that $\delta_{k}=\sum\limits_{l=1}^{n_{k+1}-1}\eta_{k+1,l}+n_{k+1}c_{1}c_{2}\cdots c_{k+1}$ and $\delta_{k}\cdot n_{1}n_{2}\cdots n_{k}=|E_{k}^{\ast}|\leq 1$. Thus
\begin{equation}
\ \frac{\log(n_{1}n_{2}\cdots n_{k})}{-\log\delta_{k}}\leq 1.\nonumber
\end{equation}
Because $\dim_{H}E=1$, using the result of Theorem \ref{Hdim}, we have the desired conclusion.
\qed
\end{prof}
Like other papers proving minimality of Cantor type sets, to estimate the Hausdorff dimension of $f(E)$, we apply the following mass distribution principle (refer to \cite{F}).
\begin{lemma}\label{Lemma-3}
Let $\mu$ be a mass distribution supported on $f(E)\subseteq \mathbb{R}^{1}$. Suppose that for some $0<d<1$, there is a number $c>0$ such that for all intervals $U\subseteq f(F_{0})$, we have $\mu(U)\leq c|U|^{d}$. Then $\dim_{H}f(E)\geq d$.
\end{lemma}
For any $0<d<1$, we are going to construct a Borel probability measure supported on $f(E)$ which satisfies the condition in Lemma \ref{Lemma-3}. Then because of the arbitrariness of $d$, we have $\dim_{H}f(E)=1$.

Let $J_{m-1}$ be a basic interval of $f(F_{m-1})$, where $m\geq 1$. Let $J_{m-1,1},\cdots,J_{m-1,N(J_{m-1})}$ be the basic intervals of $f(F_{m})\bigcap J_{m-1}$, $N(J_{m-1})$ is the number of basic intervals of $f(F_{m})\bigcap J_{m-1}$ and $N(J_{m-1})\leq4$. Set $C(J_{m-1},d)=\sum_{i=1}^{N(J_{m-1})}|J_{m-1,i}|^{d}$. Then there is a unique Borel probability measure $\mu$ supported on $f(E)$ such that for any basic interval of $f(F_{m-1})$, say $J_{m-1}$, we have
\begin{equation}
\ \mu(J_{m-1,i})=\frac{|J_{m-1,i}|^{d}}{C(J_{m-1},d)}\mu(J_{m-1}),\quad i=1,2,\cdots,N(J_{m-1}).\nonumber
\end{equation}

\begin{lemma}\label{Key-Lemma-2}
Suppose $\{b_m\}_{m\in \mathbb{N}\cup\{0\}}$ is a nonnegative real number sequence with
$$
\lim_{m\rightarrow \infty}\frac{1}{m}\sum_{i=0}^{m-1}b_i=0
$$
and  $\varepsilon\in (0,1)$. Then we have
$$
\lim_{m\rightarrow\infty}\frac{S(\varepsilon,m)}{m}=1.
$$
Here $S(\varepsilon,m)=\#(\{0\leq i\leq m-1:0\leq b_i<\varepsilon\})$.
\end{lemma}
\begin{prof}
Because $\lim\limits_{m\rightarrow\infty}\frac{1}{m}\sum_{j=0}^{m-1}b_{j}=0$, we have
\begin{align}
 \lim\limits_{m\rightarrow\infty}\frac{S(\varepsilon,m)}{m}
 &=1-\lim\limits_{m\rightarrow\infty}\frac{m-S(\varepsilon,,m)}{m}\nonumber\\
 &\geq 1-\lim\limits_{m\rightarrow\infty}\frac{1}{m\varepsilon}\sum_{j=0}^{m-1}b_{j}=1.\nonumber
\end{align}
\qed
\end{prof}

\textbf{Proof of theorem \ref{theorem-2}:} To complete the proof, we only need to show that the Borel probability measure $\mu$ constructed above satisfying the condition in Lemma \ref{Lemma-3}, i.e., there exists a constant $c$ such that for any interval $U\subseteq f(F_{0})$, we have $\mu(U)\leq c|U|^{d}$.
\\
\indent \textbf{Case 1}: First, suppose that $U$ is some basic interval of $f(F_{m})$, say $J_{m}$. For any $0\leq j\leq m-1$, let $J_{j}$ be the basic interval of $f(F_{j})$ such that
\begin{equation}
\ J_{m}\subseteq J_{m-1}\subseteq\cdots\subseteq J_{1}\subseteq J_{0}=f(F_{0})\nonumber.
\end{equation}
By the definition of $\mu$, we have
\begin{equation}
\ \frac{\mu(J_{m})}{|J_{m}|^{d}}|J_{0}|^{d}=\prod_{j=0}^{m-1}\frac{|J_{j}|^{d}}{C(J_{j},d)} .\nonumber
\end{equation}
So it suffices to show
\begin{equation}
\ \liminf\limits_{m\rightarrow\infty}(\prod_{j=0}^{m-1}\frac{C(J_{j},d)}{|J_{j}|^{d}})^{\frac{1}{m}}>1. \nonumber
\end{equation}
To accomplish this goal, we are going to estimate $\frac{C(J_{j},d)}{|J_{j}|^{d}}$ for $0\leq j\leq m-1$.
\\
\indent Let $I_{j}$ be the basic interval of $F_{j}$ such that $f(I_{j})=J_{j}$. Let $J_{j,1},\cdots,J_{j,N(J_{j})}$ be basic intervals of $f(F_{j+1})\bigcap J_{j}$ and the gaps between them in $J_{j}$ are $L_{j,0},\cdots,L_{j,N(J_{j})}$. It's worth noting that the gap may be empty.  $L_{j,0},J_{j,1},L_{j,1},\cdots,J_{j,N(J_{j})},L_{j,N(J_{j})}$ locate from left to right and their union is $J_{j}$. Let $I_{j,l}$ be the basic interval of $F_{j+1}$ such that $f(I_{j,l})=J_{j,l}$ for $1\leq l\leq N(J_{j})$. Let $G_{j,l}$ be component interval of $I_{j}-F_{j+1}$ such that $f(G_{j,l})=L_{j,l}$ for $0\leq l\leq N(J_{j})$.
\\
\indent Recognize that
\begin{equation}
\ \frac{C(J_{j},d)}{|J_{j}|^{d}}=\frac{C(J_{j},d)}{C(J_{j},1)^{d}}\frac{C(J_{j},1)^{d}}{|J_{j}|^{d}}.\nonumber
\end{equation}

Fix a sufficient small $\varepsilon>0$, actually we only need $\varepsilon$ to satisfy the following three conditions:

$(1)$ $\varepsilon<\frac{1-\alpha}{5}$;

$(2)$ $(1-20x^{p})\geq (1-x^{p})^{21}>0$ for any $x\in[0,\varepsilon)$;

$(3)$ $\log(1-x^{p})\geq -2x^{p}$ for any  $x\in[0,\varepsilon)$.

Without loss of generality, let $J_{j,1}$ be the biggest intervals among $J_{j,1},J_{j,2},\cdots,J_{j,N(J_{j})}$ and put $x_{l}=\frac{|J_{j,l}|}{|J_{j,1}|}$, then $x_{1}=1,\,0<x_{l}\leq 1$. We have
\begin{align} \frac{C(J_{j},d)}{C(J_{j},1)^{d}}&=\frac{1+x_{2}^{d}+x_{3}^{d}+\cdots+x_{N(J_{j})}^{d}}{(1+x_{2}+x_{3}+\cdots+x_{N(J_{j})})^{d}}\nonumber\\
&\geq (1+x_{2}+x_{3}+\cdots+x_{N(J_{j})})^{1-d}\nonumber\\
&\geq 1.\nonumber
\end{align}
Thus
\begin{equation}
\ \frac{C(J_{j},d)}{|J_{j}|^{d}}
=\frac{C(J_{j},d)}{C(J_{j},1)^{d}}\frac{C(J_{j},1)^{d}}{|J_{j}|^{d}}
\geq\frac{C(J_{j},1)^{d}}{|J_{j}|^{d}}.\nonumber
\end{equation}

If $\beta_{j}<\varepsilon$, by the definition of $\beta_{j}$, we have $\frac{|G_{j,l}|}{|I_{j}|}\leq\beta_{j}$. Then by Lemma \ref{Lemma-1}, we have $\frac{|L_{j,l}|}{|J_{j}|}\leq 4\beta_{j}^{p}$ for $0\leq l\leq N(J_{j})$. So
\begin{equation}\label{Mestimation}
\ (\frac{\sum\limits_{l=1}^{N(J_{j})}|J_{j,l}|}{|J_{j}|})^{d}\geq(1-20\beta_{j}^{p})^{d}\geq (1-\beta_{j}^{p})^{21d}.
\end{equation}

If $\beta_{j}<\varepsilon$ and $\gamma_{j+1}<\alpha$, by Lemma \ref{Lemma-1} and Jensen's inequality, we have
\begin{align}
 \frac{|J_{j,2}|+|J_{j,3}|+\cdots+|J_{j,N(J_{j})}|}{|J_{j}|}
 &\geq\lambda\frac{|I_{j,2}|^{q}+|I_{j,3}|^{q}+\cdots+|I_{j,N(J_{j})}|^{q}}{|I_j|^{q}}\nonumber\\
 &\geq3^{1-q}\lambda(\frac{|I_{j,2}|+|I_{j,3}|+\cdots+|I_{j,N(J_{j})}|}{|I_j|})^{q}.\nonumber
\end{align}
Since $G_{j,l}\subset I_j,\frac{|G_{j,l}|}{|I_j|}\leq\beta_j< \varepsilon$ for all $0\leq l\leq N(J_i)\leq 4$ and $\gamma_{j+1}<\alpha$, we have
\begin{align}
\frac{|I_{j,2}|+|I_{j,3}|+\cdots+|I_{j,N(J_{j})}|}{|I_j|}
&=\frac{|I_j|-|I_{j,1}|-\sum_{l=0}^{N(J_j)}|G_{j,l}|}{|I_j|}\nonumber\\
&\geq 1-\alpha-5\varepsilon.\nonumber
\end{align}
From the above two inequalities, we have
$$
\frac{|J_{j,2}|+|J_{j,3}|+\cdots+|J_{j,N(J_{j})}|}{|J_{j}|}\geq3^{1-q}\lambda(1-5\varepsilon-\alpha)^{q}.
$$

By Lemma \ref{Lemma-1}, we obtain
$$
\frac{|J_{j,1}|}{|J_j|}=\frac{|f(I_{j,1})|}{|f(I_j)|}\leq4\frac{|I_{j,1}|^p}{|I_j|^p}\leq 4\alpha^{p},
$$
then we have
\begin{equation}
\ x_{2}+x_{3}+\cdots+x_{N(J_{j})}\geq\frac{|J_j|}{|J_{j,1}|}\frac{\lambda(1-5\varepsilon-\alpha)^{q}}{3^{q-1}}
\geq\frac{\lambda(1-5\varepsilon-\alpha)^{q}}{4\alpha^{p}\cdot3^{q-1}}.\nonumber
\end{equation}

In conclusion,
if $\beta_{j}<\varepsilon$, we have
\begin{equation}\label{eq-add}
\ \frac{C(J_{j},d)}{|J_{j}|^{d}}\geq (1-\beta_{j}^{p})^{21d};
\end{equation}
if $\beta_{j}<\varepsilon, \gamma_{j+1}<\alpha$, we have
\begin{equation}\label{eq-1}
\ \frac{C(J_{j},d)}{|J_{j}|^{d}}\geq \eta(1-\beta_{j}^{p})^{21d},
\end{equation}
where $\eta=(1+\frac{\lambda(1-5\varepsilon-\alpha)^{q}}{4\alpha^{p}\cdot3^{q-1}})^{1-d}>1$.

Besides, if $\beta_{j}<\varepsilon$,
\begin{align}
\ 0\geq\frac{1}{m}\sum_{\substack{j=0\\ \beta_{j}<\varepsilon}}^{m-1}  \log(1-\beta_{j}^{p})
&\geq \frac{-2}{m}\sum_{\substack{j=0\\   \beta_{j}<\varepsilon}}^{m-1}\beta_{j}^{p}\geq \frac{-2}{m}\sum_{j=0}^{m-1}\beta_{j}^{p}\nonumber \\
&\geq -2(\frac{1}{m}\sum_{j=0}^{m-1}\beta_{j})^{p}\rightarrow 0,\quad\text{as $m\rightarrow\infty$}.\nonumber
\end{align}
So we obtain
\begin{equation}\label{eq-3}
\lim\limits_{m\rightarrow\infty}[\prod_{\substack{j=0\\   \beta_{j}<\varepsilon}}^{m-1}(1-\beta_{j}^{p})]^{\frac{1}{m}}=1.
\end{equation}

If $\beta_j\geq\varepsilon$, by Lemma \ref{Lemma-1}, we have
\begin{align}
\frac{C(J_{j},1)}{|J_{j}|}=\frac{\sum\limits_{l=1}^{N(J_{j})}|J_{j,l}|}{|J_{j}|}
 &\geq \lambda \frac{\sum\limits_{l=1}^{N(J_{j})}|I_{j,l}|^{q}}{|I_{j}|^{q}}\geq \frac{\lambda}{4^{q-1}}(\frac{\sum\limits_{l=1}^{N(J_{j})}|I_{j,l}|}{|I_{j}|})^{q}\nonumber\\
 &\geq \frac{\lambda}{4^{q-1}}\Gamma_j^{q}.\nonumber
\end{align}
So we obtain
\begin{equation}\label{eq-2}
\ \frac{C(J_{j},d)}{|J_{j}|^{d}}\geq\frac{C(J_{j},1)^{d}}{|J_{j}|^{d}}\geq (\frac{\lambda}{4^{q-1}}\Gamma_j^{q})^{d}.
\end{equation}

For every $m\geq1$, let $S(m)=\#(\{0\leq j\leq m-1:\beta_{j}<\varepsilon\})$, $T(m)=\#(\{1\leq j\leq m:\gamma_{j}<\alpha\})$ and $ST(m)=\#(\{1\leq j\leq m:\beta_{j-1}<\varepsilon,\,\gamma_{j}<\alpha\})$. Because $\lim\limits_{m\rightarrow\infty}\frac{1}{m}\sum_{j=0}^{m-1}\beta_{j}=0$, by Lemma \ref{Key-Lemma-2} we get
\begin{equation}
 \lim\limits_{m\rightarrow\infty}\frac{S(m)}{m}=1.\nonumber
\end{equation}
On the one hand, we assume that $\liminf\limits_{m\rightarrow\infty}\frac{T(m)}{m}=t>0$, then we get
\begin{equation}\label{eq-4}
\ \liminf\limits_{m\rightarrow\infty}\frac{ST(m)}{m}\geq t.
\end{equation}
Combining equation (\ref{eq-add}), (\ref{eq-1}), (\ref{eq-3}), (\ref{eq-2}) and (\ref{eq-4}), we obtain
\begin{equation*}
\begin{split}
\prod_{j=0}^{m-1}\frac{C(J_{j},d)}{|J_{j}|^{d}}
&=\prod_{\substack{j=0\\   \beta_{j}<\varepsilon,\gamma_{j+1}<\alpha}}^{m-1}\frac{C(J_{j},d)}{|J_{j}|^{d}}
\prod_{\substack{j=0\\ \beta_{j}<\varepsilon,\gamma_{j+1}\geq\alpha}}^{m-1}\frac{C(J_{j},d)}{|J_{j}|^{d}}
\prod_{\substack{j=0\\ \beta_{j}\geq\varepsilon,}}^{m-1}\frac{C(J_{j},d)}{|J_{j}|^{d}}\\
&\geq \eta^{ST(m)}\prod_{\substack{j=0\\ \beta_{j}<\varepsilon}}^{m-1}(1-\beta_{j}^{p})^{21d}
\prod_{\substack{j=0\\ \beta_{j}\geq\varepsilon,}}^{m-1}(\frac{\lambda}{4^{q-1}}\Gamma_j^{q})^{d}\\
&\geq\eta^{ST(m)}\prod_{\substack{j=0\\ \beta_{j}<\varepsilon}}^{m-1}(1-\beta_{j}^{p})^{21d}
(\prod_{j=0}^{m-1}\Gamma_j)^{qd}
\prod_{\substack{j=0\\ \beta_{j}\geq\varepsilon,}}^{m-1}(\frac{\lambda}{4^{q-1}})^{d}\\
&=\eta^{ST(m)}\prod_{\substack{j=0\\ \beta_{j}<\varepsilon}}^{m-1}(1-\beta_{j}^{p})^{21d}
(\prod_{j=0}^{m-1}\Gamma_j)^{qd}
(\frac{\lambda}{4^{q-1}})^{d(m-S(m))}.\\
\end{split}
\end{equation*}
From the above inequality and equation (\ref{eq-3}), we have
\begin{align}
\ \liminf\limits_{m\rightarrow\infty}(\prod_{j=0}^{m-1}\frac{C(J_{j},d)}{|J_{j}|^{d}})^{\frac{1}{m}}
&\geq \liminf\limits_{m\rightarrow\infty}\eta^{\frac{ST(m)}{m}}\lim\limits_{m\rightarrow\infty}(\prod_{j=0   }^{m-1}\Gamma_j)^{\frac{dq}{m}}(\frac{\lambda}{4^{q-1}})^{\frac{d(m-S(m))}{m}} \nonumber\\
&> 1+r,\nonumber
\end{align}
where $1<r+1<\eta^{t}$. To conclude, we have found some constant $C>0$ such that for any $m\geq0$ and any $I_m\in\mathcal{F}_{m}$, $J_m=f(I_m)$,
\begin{equation}\label{eq-5}
\ \mu(J_m)\leq C\frac{|J_m|^{d}}{(1+r)^{m}}.
\end{equation}
\\
\indent \textbf{Case 2}: After case 1, we assume that $J$ is an interval such that $J\subseteq f(F_0)$. Without losing generality, we assume that $f(F_0)\not\subseteq J$. Let $\Phi_k^{J}=\{I\in\mathcal{F}_{k}:f(I)\subseteq J,\,f(Fa(I))\not\subseteq J\}$, recall step 2 of the reconstruction, we have  $\#(\Phi_k^{J})\leq 8$ for any $k\geq 1$.
Since $0<d<1$, by Jensen's inequality, we have
$$
\frac{1}{\#(\Phi_{k}^{J})}\Sigma_{I\in\Phi_{k}^{J}}|f(I)|^{d}\leq (\frac{\Sigma_{I\in\Phi_{k}^{J}}|f(I)|}{\#(\Phi_{k}^{J})})^{d}.
$$
According to the definition of $\Phi_{k}^{J}$, we have
$$
(\Sigma_{I\in\Phi_{k}^{J}}|f(I)|)^d\leq |J|^d.
$$
Thus
\begin{align}
 \mu(J)
 &\leq\mu(\cup_{k\geq1}\cup_{I\in\Phi_k^{J}}f(I))=\Sigma_{k\geq1}\Sigma_{I\in\Phi_{k}^{J}}\mu(f(I))\nonumber\\
 &\leq C\Sigma_{k\geq1}\Sigma_{I\in\Phi_{k}^{J}}\frac{|f(I)|^{d}}{(1+r)^{k}}\leq 8^{1-d}C\Sigma_{k\geq1}\frac{(\Sigma_{I\in\Phi_{k}^{J}}|f(I)|)^{d}}{(1+r)^{k}}\nonumber\\
&\leq 8^{1-d}C \Sigma_{k\geq1}\frac{|J|^d}{(1+r)^{k}}= \frac{8^{1-d}C}{r}|J|^d.\nonumber
\end{align}
Since $d$ is  arbitrary, by Lemma \ref{Lemma-3}, we have $\dim_{H}f(E)=1$,
which completed the proof.
\qed

When studying the paper \cite{WW}, we have found some easier way to check whether homogeneous perfect sets are minimal or not than Theorem \ref{theorem-2}, see Theorem 5.1 in \cite{WW}. We will show below that it is just the special case of Theorem \ref{theorem-2}. Before introducing our result, we define some notations. For any $m\geq 1$, set
$$\Lambda_{m}=\frac{\max_{I\in\mathcal{F}_{m}}|I|}{\min_{I\in\mathcal{F}_{m-1}}|I|},\quad \lambda_{m}=\frac{\min_{I\in\mathcal{F}_{m}}|I|}{\max_{I\in\mathcal{F}_{m-1}}|I|}.
$$

\begin{theorem}\label{theorem-3}
Let $E=\bigcap_{m\geq0}F_{m}$ be a homogeneous perfect set, then $E$ is minimal if it satisfies the conditions $(a)$, $(b)$, $(c)$:

$(a)$ $\lim\limits_{m\rightarrow\infty}\frac{1}{m}\log_{2}|F_{m}|=0$

$(b)$ $\lim\limits_{m\rightarrow\infty}\frac{1}{m}\sum_{j=0}^{m-1}\beta_{j}=0$

$(c)$ there exists $K\geq 1$ such that $K \min_{I\in\mathcal{F}_{m}}|I|\geq \max_{I\in\mathcal{F}_{m}}|I|$ for all $m\geq 1$.
\end{theorem}
\begin{remark}We prove this theorem by Theorem \ref{theorem-2}. But in the last section, we introduce an example in which the homogeneous perfect set is minimal while some condition in Theorem \ref{theorem-3} does not hold.
\end{remark}
\begin{prof}
Obviously, the set $E$ satisfies condition $(\tilde{a})$ in Theorem \ref{theorem-2}. It suffices to show that $E$ satisfies the conditions $(\tilde{b})$ and $(\tilde{c})$ in Theorem \ref{theorem-2}. Fix $\varepsilon\in(0,\frac{1}{5})$ such that $\log(1-5x)\geq -10x$ for any  $x\in[0,\varepsilon)$, and let $S(m)=\#(\{0\leq j\leq m-1:\beta_{j}<\varepsilon\})$.

For $\beta_{j}<\varepsilon$, similar to the proof in Theorem \ref{theorem-2}, we have
\begin{equation}\label{eq-6}
\ \lim\limits_{m\rightarrow\infty}(\prod_{\substack{j=0\\   \beta_{j}<\varepsilon}}^{m-1}(1-5\beta_j))^{\frac{1}{m}}=1.\nonumber
\end{equation}

According to condition $(c)$, there exists $K\geq 1$ such that
$$
K \min_{I\in\mathcal{F}_{m}}|I|\geq \max_{I\in\mathcal{F}_{m}}|I|
$$ for all $m\geq 1$. Thus we have $K^{2}\lambda_{j+1}\geq \Lambda_{j+1}$.

Let $J\in \mathcal{F}_{j}$ such that $\Gamma_j=\frac{C(J,1)}{|J|}$. It is easy to see that
$$
\Gamma_j=\frac{C(J,1)}{|J|}\geq \frac{\min_{I\in \mathcal{F}_{j+1}} |I|}{\max_{I\in \mathcal{F}_j}|I|}=\lambda_{j+1}.
$$
On the other hand, because every basic interval of $F_{j-1}$ has at most $4$ basic intervals in $F_{j}$, we have $\frac{|F_{j}|}{|F_{j-1}|}\leq\min\{1,4\Lambda_{j}\}$ for any $1\leq j\leq m$. Then for any subset $\Delta\subseteq\{1,2,\cdots,m\}$, we have $|F_{m}|\leq\prod_{j\in\Delta}(4\Lambda_{j})$. Thus for any subset $\Delta\subseteq\{0,1,\cdots,m-1\}$, we have
\begin{equation}\label{eq-6}
\ \prod_{j\in\Delta}(4K^{2}\Gamma_{j})\geq |F_{m}|.\nonumber
\end{equation}

Still by Lemma \ref{Key-Lemma-2}, we have
\begin{equation*}
\lim_{m\rightarrow\infty}(1-\frac{S(m)}{m})=0.
\end{equation*}
By Lemma \ref{key-lemma}, we have $\Gamma_j\geq 1-5\beta_{j}$.
Then we get
\begin{align}
 \lim\limits_{m\rightarrow\infty}(\prod_{j=0}^{m-1}\Gamma_j)^{\frac{1}{m}}
 &=\lim\limits_{m\rightarrow\infty}(\prod_{\substack{j=0\\   \beta_{j}<\varepsilon}}^{m-1}\Gamma_j)^{\frac{1}{m}}(\prod_{\substack{j=0\\   \beta_{j}\geq\varepsilon}}^{m-1}\Gamma_j)^{\frac{1}{m}}
 \nonumber\\
 &\geq\lim\limits_{m\rightarrow\infty}(\prod_{\substack{j=0\\   \beta_{j}<\varepsilon}}^{m-1}(1-5\beta_j))^{\frac{1}{m}}(\prod_{\substack{j=0\\   \beta_{j}\geq\varepsilon}}^{m-1}\frac{1}{4K^{2}})^{\frac{1}{m}}|F_{m}|^{\frac{1}{m}}\nonumber\\
 &=\lim\limits_{m\rightarrow\infty}(\prod_{\substack{j=0\\   \beta_{j}<\varepsilon}}^{m-1}(1-5\beta_j))^{\frac{1}{m}}(\frac{1}{4K^2})^{1-\frac{S(m)}{m}}|F_m|^{\frac{1}{m}}\nonumber\\
 &=1.\nonumber
\end{align}
That is $E$ satisfies the condition $(\tilde{b})$ in Theorem \ref{theorem-2}.

By condition $(c)$, for any $J\in\mathcal{F}_{j}$, we have
$$
\gamma_j\leq \frac{\max_{J\in \mathcal{F}_j}|J|}{\min_{\widetilde{J}\in \mathcal{F}_{j-1}}|\widetilde{J}|}\leq \frac{K|J|}{\frac{1}{K}\max_{\widetilde{J}\in \mathcal{F}_{j-1}}|\widetilde{J}|}\leq K^{2}\frac{|J|}{|Fa(J)|}.
$$
Suppose that $\gamma_j=\frac{|J|}{|Fa(J)|}$ for some $J\in \mathcal{F}_{j}$.
Since $Fa(J)$ contains at least $2$ elements in $\mathcal{F}_{j}$, we have
$$
1\geq \gamma_j+\frac{|J^{*}|}{|Fa(J)|}=\gamma_j+\frac{|J^{*}|}{|Fa(J^{*})|}\geq\gamma_j+\frac{\gamma_j}{K^2}
$$
for some $J^{*}\subset Fa(J)$.
Take $\alpha\in(\frac{K^{2}}{K^{2}+1} ,1)$, we have
$$
\liminf\limits_{m\rightarrow\infty}\frac{\#(\{1\leq i\leq m:\,\gamma_i<\alpha\})}{m}=1>0.
$$\qed
\end{prof}

\section{Proof of Theorem \ref{theorem-1}}
Now, we prove Theorem \ref{theorem-1} in this section.
According to Theorem \ref{theorem-3}, we need to prove the conditions $(a)$, $(b)$, $(c)$ hold in Case $(A)$ and Case $(B)$. Before our checking, we introduce some notations for concision. For any $k\in \mathbb{N}$,
$$
B(k)\equiv\max_{1\leq l\leq n_{k}-1}\eta_{k,l}+\eta_{k+1,0}+\eta_{k+1,n_{k+1}}, A(k)\equiv\max_{1\leq l\leq n_{k}-1}\eta_{k,l}
$$
and
$$
b(k)\equiv\min_{1\leq l\leq n_{k}-1}\eta_{k,l}+\eta_{k+1,0}+\eta_{k+1,n_{k+1}}, \eta(k)\equiv\eta_{k+1,0}+\eta_{k+1,n_{k+1}}.
$$
For any $I\in\mathcal{F}_{m}$, let $m_{k-1}\leq m<m_{k}$ for some $k\geq1$, we set
$$
\chi(I,\mathcal{F}_{m_{k}})=\#(\{J\in \mathcal{F}_{m_{k}}: J\subset I\}).
$$
That is $\chi(I,\mathcal{F}_{m_{k}})$ is the number of basic intervals of $F_{m_{k}}$ that are contained in $I$. Finally, let $I(m),I(M)$ be the members of $\mathcal{F}_m$ such that
$$
\max_{I\in \mathcal{F}_m}{|I|}=|I(M)|\quad\mathrm{and}\quad\min_{I\in \mathcal{F}_m}{|I|}=|I(m)|.
$$
\\
\indent \textbf{Case (A)}.
Condition (a): First, if $m=m_{k}$ for some $k\geq1$, then we have
$$
|F_m|=|E_{k}^{*}|=n_{1}n_{2}\cdots n_{k}\delta_{k}.
$$
Thus,
\begin{align}
\ \frac{1}{m_{k}}\log_{2}(n_{1}n_{2}\cdots n_{k}\delta_{k})
& =\frac{\log_{2}(n_{1}n_{2}\cdots n_{k})}{m_{k}}\cdot\frac{\log_{2}(n_{1}n_{2}\cdots n_{k})+\log_{2}\delta_{k}}{\log_{2}(n_{1}n_{2}\cdots n_{k})}\nonumber\\
& \geq2(1-(\frac{\log_{2}(n_{1}n_{2}\cdots n_{k})}{-\log_{2}\delta_{k}})^{-1}) \rightarrow 0\nonumber
\end{align}
as $k\rightarrow\infty$. Here we use the fact that $\log_{2}n_{j}\leq i_{j}+1$ for $1\leq j\leq k$ and $\dim_{H}E=1$.
\\
\indent Now suppose that $m_{k-1}<m<m_{k}$, we consider $F_{m_{k-1}+i_{k}-1}$. Actually for every $\sigma\in\Omega_{k-1}$, the union of all basic intervals of $F_{m_{k-1}+i_{k}-1}\bigcap J_{\sigma}^{\ast}$ may be seen as $J_{\sigma}^{\ast}$ subtracting $2$ edging gaps of $J_{\sigma}^{\ast}$, that is the gaps of lengths $\eta_{k,0}^{\ast}$, $\eta_{k,n_{k}}^{\ast}$ and subtracting $2^{i_{k}-1}-1$ middle gaps of $J_{\sigma}^{\ast}$. The gaps of $J_{\sigma}^{\ast}$ mean the component intervals of $J_{\sigma}^{\ast}-(J_{\sigma}^{\ast}\bigcap F_{m_k})$. Then the union still obtain at least $2^{i_{k}-1}$ middle gaps of $J_{\sigma}^{\ast}$. Each remaining middle gap has length at least $b(k)$ and each subtracted gap has length at most $B(k)$. With the equation $Cb(k)\geq B(k)$ and $\eta_{k,0}^{\ast}+\eta_{k,n_{k}}^{\ast}\leq b(k)$, we have
\begin{equation}
\ |F_{m_{k-1}+i_{k}-1}|(1+C)\geq |F_{m_{k-1}}|.\nonumber
\end{equation}
Thus for any $m_{k-1}<m<m_{k}$, we have $|F_{m}|(1+C)\geq |F_{m_{k-1}}|$.

For any $\varepsilon>0$, there exists $N>0$ such that $\frac{1}{m_{k}}\log_{2}|F_{m_{k}}|>-\frac{\varepsilon}{2}$ for any $k\geq N$ and $\frac{1}{N}\log_{2}(1+C)<\frac{\varepsilon}{2}$. So if $m>m_{N}$, there exists a unique $r\geq N$ such that $m_{r}\leq m<m_{r+1}$, thus $\frac{1}{m}\log_{2}|F_{m}|\geq \frac{1}{m_{r}}(\log_{2}|F_{m_{r}}|+\log_{2}(1+C)^{-1})>-\varepsilon$. In conclude, condition (a) is satisfied.
\\
\indent Condition (c): If $m=m_k$ for some $k\in \mathbb{N}$, we have
$$
\min_{I\in \mathcal{F}_m}|I|=\max_{I\in \mathcal{F}_m}|I|=\delta_k.
$$
Thus condition $(c)$ holds in this case.
\\
\indent Now let $m_{k-1}< m< m_{k}$ for some $k\geq1$. Each basic interval of $F_{m}$ may been seen as the union of some basic intervals of $F_{m_{k}}$ and the gaps between these basic intervals. While by the reconstruction of homogeneous perfect set, the numbers of the basic intervals of $F_{m_{k}}$ which are contained in them differ at most one, that is we have $\chi(I(M),\mathcal{F}_{m_{k}})\leq\chi(I(m),\mathcal{F}_{m_{k}})+1$. Thus
\begin{align}
\ \max_{I\in\mathcal{F}_{m}}|I|=|I(M)|
& \leq (\chi(I(m),\mathcal{F}_{m_{k}})+1)\delta_{k}+\chi(I(m),\mathcal{F}_{m_{k}})B(k)\nonumber\\
& \leq 2C[\chi(I(m),\mathcal{F}_{m_{k}})\delta_{k}+(\chi(I(m),\mathcal{F}_{m_{k}})-1)b(k)]\nonumber\\
& \leq 2C|I(m)|=2C\min_{I\in\mathcal{F}_{m}}|I|.\nonumber
\end{align}
Thus condition $(c)$ holds for $K=2C$.
\\
\indent Condition (b): For any $m_{k-1}\leq m<m_{k}$, denote $\tilde{\gamma}_{m}=\min\{\frac{\eta_{k,l}^{\ast}}{|I|}: I\in\mathcal{F}_{m}, 1\leq l\leq n_{k}-1\}$. By condition $(c)$ we have proved above and condition $(A)$ of Theorem \ref{theorem-1}, it is easy to see $2C^{2}\tilde{\gamma}_{m}\geq\beta_{m}$. Besides, for any $0\leq j\leq m-1$ and any basic interval $I$ of $F_{j}$, to get $I\cap F_{j+1}$, $I$ must subtract a interval of length at least $\tilde{\gamma}_{j}|I|$, that is $\frac{|F_{j+1}|}{|F_{j}|}\leq 1-\tilde{\gamma}_{j}$. Through this we get
\begin{equation}
\ |F_{m}|\leq |J_{\emptyset}^{\ast}|\prod_{j=0}^{m-1}(1-\tilde{\gamma}_{j}).\nonumber
\end{equation}
Then
\begin{equation}
\ 0\geq -\frac{1}{m}\sum\limits_{j=0}^{m-1}\tilde{\gamma}_{j}\geq\frac{1}{m}\sum\limits_{j=0}^{m-1}\log(1-\tilde{\gamma}_{j})\rightarrow 0\nonumber
\end{equation}
as $m\rightarrow\infty$.
We have completed the proof.

\textbf{Case (B)}.
Condition (a). As we did in Case $(A)$, we have
\begin{align}
\lim\limits_{k\rightarrow\infty}\frac{1}{m_{k}}\log_{2}|F_{m_{k}}|
&=\lim\limits_{k\rightarrow\infty}\frac{1}{m_{k}}\log_{2}(n_1n_2\cdots n_k\delta_k)=0.\nonumber
\end{align}
Now suppose $m_{k-1}<m<m_{k}$, we consider $F_{m_{k-1}+i_{k}-1}$. We try to compare $|F_{m_{k-1}+i_{k}-1}|$ with $|F_{m_{k-1}}|$. For any $\sigma\in\Omega_{k-1}$, we compare $J_{\sigma}^{\ast}\cap F_{m_{k-1}+i_{k}-1}$ and $J_{\sigma}^{\ast}$. Very similar to condition $(a)$ in case $(A)$, $J_{\sigma}^{\ast}\cap F_{m_{k-1}+i_{k}-1}$ contains at least $2^{i_{k}-1}$ basic intervals of $F_{m_{k}}$ of length $\delta_{k}$ and contains at least $2^{i_{k}-1}$ middle gaps of $J_{\sigma}^{\ast}$ with length at least $\eta(k)$. Recall that gaps of $J_{\sigma}^{\ast}$ denote the component intervals of $J_{\sigma}^{\ast}-(J_{\sigma}^{\ast}\cap F_{m_k})$. It is worth noting that $\eta(k)$ may be zero. On the other hand, $J_{\sigma}^{\ast}$ subtracts at most $2^{i_{k}-1}-1$ middle gaps of $J_{\sigma}^{\ast}$ with length at most $B(k)$ and two edging gaps of lengths $\eta_{k+1,0}$ and $\eta_{k+1,n_{k+1}}$, and then we get $J_{\sigma}^{\ast}\cap F_{m_{k-1}+i_{k}-1}$. Besides, $c_1c_2\cdots c_k=\delta_k+\eta_{k+1,0}+\eta_{k+1,n_{k+1}}$ and $\max_{1\leq l\leq n_{k}-1}\eta_{k,l}\leq D c_{1}c_{2}\cdots c_{k}$, thus we get $A(k)\leq D(\delta_k+\eta(k))$. But
\begin{align}
(2^{i_{k}-1}-1)B(k)+\eta(k)
&\leq2^{i_{k}-1}(A(k)+\eta(k))\nonumber\\
&\leq2^{i_{k}-1}(1+D)(\delta_{k}+\eta(k)),\nonumber
\end{align}
we have $|F_{m_{k-1}+i_{k}-1}|(2+D)\geq|F_{m_{k-1}}|$. Using the same argument of case $(A)$, we have proved condition $(a)$ is satisfied.

Condition (c): We may suppose $m_{k-1}< m< m_{k}$ for some $k\in \mathbb{N}$, while the proof is trivial when $m=m_{k-1}$ according to case $(A)$. We just need the following equation
\begin{align}
\max_{I\in\mathcal{F}_{m}}|I|
&=|I(M)|\nonumber\\
&\leq\chi(I(M),\mathcal{F}_{m_{k}})\delta_k+(\chi(I(M),\mathcal{F}_{m_{k}})-1)(A(k)+\eta(k))\nonumber\\
&\leq(1+D)[(\chi(I(m),\mathcal{F}_{m_{k}})+1)\delta_k+\chi(I(m),\mathcal{F}_{m_{k}})\eta(k)]\nonumber\\
&\leq2(1+D)[\chi(I(m),\mathcal{F}_{m_{k}})\delta_k+(\chi(I(m),\mathcal{F}_{m_{k}})-1)\eta(k)]\nonumber\\
&\leq2(1+D)|I(m)|=2(1+D)\min_{I\in \mathcal{F}_m}|I|.\nonumber
\end{align}
Thus condition $(c)$ holds for $K=2(1+D)$.

Condition (b). By condition $(a)$, if $k\rightarrow\infty$, we have
\begin{align} &\frac{1}{m_{k}}\log\prod_{j=1}^{k}\frac{\delta_{j-1}-n_{j}(\eta_{j+1,0}+\eta_{j+1,n_{j+1}})-(\eta_{j,1}+\cdots+\eta_{j,n_{j}-1} )}{\delta_{j-1}}\nonumber\\
&=\frac{1}{m_{k}}\log\prod_{j=1}^{k}\frac{n_{j}\delta_{j}}{\delta_{j-1}} =\frac{1}{m_{k}}\log(n_1\cdots n_k\delta_k)-\frac{1}{m_{k}}\log\delta_0\nonumber\\
&=\frac{1}{m_{k}}\log |F_{m_k}|-\frac{1}{m_{k}}\log\delta_0\rightarrow 0.\nonumber
\end{align}

Using the fact that $\log(1-x)\leq-x$ for any $x\in[0,1)$, if $k\rightarrow\infty$,
\begin{equation}
\  \frac{1}{m_{k}}\sum_{j=1}^{k}\frac{\eta_{j,1}+\cdots+\eta_{j,n_{j}-1}+n_{j}(\eta_{j+1,0}+\eta_{j+1,n_{j+1}})}{\delta_{j-1}}\rightarrow0.\nonumber
\end{equation}
Then it is obvious that if $k\rightarrow\infty$
\begin{equation}\label{b-1}
\  \frac{1}{m_{k}}\sum_{j=1}^{k}\frac{B(j)}{\delta_{j-1}}\rightarrow0.
\end{equation}

Now, we begin to estimate $\beta_m$ for $m\geq 0$. Let $m_{k-1}\leq m<m_{k}$ for some $k\in \mathbb{N}$.
If $I\in \mathcal{F}_{m_{k}-1}$, then $I$ contains at least $2$ element in $\mathcal{F}_{m_{k}}$, thus $|I|\geq 2\delta_k+\eta(k)$. If $I\in \mathcal{F}_{m_{k}-2}$, then $I$ contains at least $2^2$ elements in $\mathcal{F}_{m_k}$, thus $|I|\geq 2^2\delta_k+(2^{2}-1)\eta(k)$. Inductively, if $s\in \{1,\cdots,m_{k}-m_{k-1}\}$, $I$ contains at least $2^{s}$ elements in $\mathcal{F}_{m_k}$ for any $I\in F_{m_k-s}$. Thus $|I|\geq 2^{s}\delta_k+(2^{s}-1)\eta(k)$. Obviously, for any $L\in \mathcal{G}_m$, we have
$$
|L|\leq B(k)\equiv\max_{1\leq l\leq n_{k}-1}\eta_{k,l}+\eta_{k+1,0}+\eta_{k+1,n_{k+1}}.
$$
So for any $s\in \{1,\cdots,m_{k}-m_{k-1}=i_k\}$,
\begin{align}
\beta_{m_k-s}&\leq \frac{B(k)}{2^{s-1}(\delta_k+\eta(k))}.
\end{align}

Thus
\begin{align}\label{c11}
\sum_{m=m_{k-1}}^{m_k-1}\beta_m&=\sum_{s=1}^{i_k}\beta_{m_k-s}\leq \frac{B(k)}{\delta_k+\eta(k)}\sum_{s=0}^{i_{k}-1}\frac{1}{2^s}\leq \frac{2B(k)}{\delta_k+\eta(k)}.
\end{align}
By the above analysis, for $k\geq1$, we have
\begin{equation}\label{Mainestimate}
\begin{split}
\frac{1}{m_{k}}\sum_{j=0}^{m_{k}-1}\beta_{j}
&\leq \frac{2}{m_{k}}\sum_{j=1}^{k}\frac{B(j)}{\delta_j+\eta(j)}.
\end{split}
\end{equation}

Now, for any $\varepsilon>0$, there exists $\delta>0$ such that
\begin{equation}\label{c-4}
0<\frac{1+D}{\log_{2}\frac{1}{(1+D)\delta}-1}<\frac{\varepsilon}{4}.
\end{equation}

Thus, if $\frac{\delta_{j}+\eta(j)}{\delta_{j-1}}<\delta$ for some $j\geq1$, using the fact that $A(j)\leq D(\delta_{j}+\eta(j))$, we have $(1+D)n_{j}\delta\geq1$. Using the fact that $i_{j}\geq\log_{2}n_{j}-1$ and equation (\ref{c-4}) we have
\begin{equation}
\  \frac{1+D}{i_{j}}<\frac{\varepsilon}{4}.
\end{equation}
By equation (\ref{b-1}), there exists $M>0$ such that for all $k\geq M$, we have
\begin{equation}\label{eq-17} \frac{1}{m_{k}}\sum_{j=1}^{k}\frac{B(j)}{\delta_{j-1}}<\frac{\varepsilon\delta}{4}.
\end{equation}
Then if $k\geq M$, we can get
\begin{align}
\
&\frac{1}{m_{k}}\sum_{j=1}^{k}\frac{B(j)}{\delta_{j}+\eta(j)} \nonumber\\
&\leq \frac{1}{m_{k}}(\sum_{\substack{j=1\\   \frac{\delta_{j}+\eta(j)}{\delta_{j-1}}<\delta}}^{k}(1+D)+\sum_{\substack{j=1\\   \frac{\delta_{j}+\eta(j)}{\delta_{j-1}}\geq\delta}}^{k}\frac{B(j)}{\delta_{j}+\eta(j)})\nonumber\\
&\leq \frac{1}{m_{k}}\sum_{j=1}^k\frac{i_{j}\varepsilon}{4}+\frac{1}{m_{k}}\sum_{j=1}^{k}\frac{B(j)}{\delta_{j-1}}\cdot\frac{1}{\delta}
\nonumber\\
&<\frac{\varepsilon}{4}+\frac{\varepsilon}{4}=\frac{\varepsilon}{2}\nonumber.
\end{align}
Then we get if $k\rightarrow\infty$, $\frac{1}{m_{k}}\sum_{j=1}^{k}\frac{B(j)}{\delta_{j}+\eta(j)}\rightarrow0$. Thus we have
\begin{align}\label{mainest}
\frac{1}{m_{k}}\sum_{j=0}^{m_{k}-1}\beta_{j}\rightarrow 0
\end{align}
as $k\rightarrow\infty$.

If $m_{k-1}<m<m_{k}$, it is easy to see that
\begin{align}\label{c-5}
\frac{1}{m}\sum_{j=0}^{m-1}\beta_{j}&= \frac{1}{m}(\sum_{j=0}^{m_{k-1}-1}\beta_{j}+\sum_{j=m_{k-1}}^{m-1}\beta_{j})\\
&\leq\frac{1}{m_{k-1}}\sum_{j=0}^{m_{k-1}-1}\beta_{j}+\frac{1}{m}\sum_{j=m_{k-1}}^{m-1}\beta_{j}\nonumber.
\end{align}
But
\begin{align}
\sum_{j=m_{k-1}}^{m-1}\beta_{j}\leq \frac{2B(k)}{\delta_k+\eta(k)}\leq 2(1+D).
\end{align}
We have the desired result
$$
\lim\limits_{m\rightarrow\infty}\frac{1}{m}\sum_{j=0}^{m-1}\beta_{j}=0,
$$
that is conditions $(b)$ holds in this case.

\section{An example}
In the end, we construct a homogeneous perfect set which satisfies the condition $(3)$ of Theorem \ref{Hdim} and does not satisfy the condition $(A)$ or $(B)$ in Theorem \ref{theorem-1}, but it is still minimal. Actually, we are aware of that this set does not satisfy some condition in Theorem \ref{theorem-3} as well. The notations in this section will be the same as notations in Section $2$.

Let $E=E(J_{\emptyset},\{n_{k}\},\{c_{k}\},\{\eta_{k,j}\})$ be a homogeneous perfect set. For any $k\geq 1$, take $c_{k}=(2^{k+1}+k(\sqrt{2})^{k})^{-1}$
and $n_{k}=2^{k}$, let $\eta_{k,0}=\eta_{k,n_k}=0$, that is the leftmost and rightmost gaps in every basic intervals are empty. Then we have
\begin{equation}
\delta_k=c_1c_2\cdots c_k.\nonumber
\end{equation}
For $k=1$, let $\eta_{k,1}=1-2\delta_k$. For $k\geq2$, set
\begin{equation}
\  l_k=\sum_{i=1}^{[\frac{k}{2}]+1}2^{k-i}=2^k-2^{k-[\frac{k}{2}]-1},\nonumber
\end{equation}
where $[\frac{k}{2}]$ is the integer part of $\frac{k}{2}$. If $l\in\{1, 2, \cdots, n_k-1\}-\{l_k\}$, let $\eta_{k,l}=\delta_k$, if $l=l_k$, let $\eta_{k,l}=2\delta_k+k(\sqrt{2})^{k}\delta_k$. We have completed the construction of the homogeneous perfect set, denote it by $E=\cap_{k\geq0}E_k$. After reconstruction in Section $2$, we have $E=\cap_{m\geq0}F_m$, where $F_{m_k}=E_k$. As in Section $3$, let $\mathcal{F}_{m}$ denote the family of all basic intervals of $F_{m}$. Recall the reconstruction in Section $2$, we have $i_{k}=k$ for $k\geq1$ and $m_k=1+2+\cdots+k$.

By the construction of $E$, for every $I\in \mathcal{F}_{m_{k-1}}$, $I$ contains exactly $n_k=2^k$ basic intervals which belong to $\mathcal{F}_{m_{k}}$ and $n_k-1$ gaps between them. We denote these gaps by $\{G_1\cdots G_{n_k-1}\}$ and $G_{l}$ lies in the left side of $G_{l+1}$ for all $l$. Then by the definition, we have
\begin{equation*}
|G_l|=\left\{
       \begin{array}{ll}
       \delta_k & \text{ if $l\in\{1,\cdots,n_k-1\}-\{l_k\}$}, \\
       2\delta_k+k(\sqrt{2})^{k}\delta_k& \text{ if $l=l_k$}.
        \end{array}
     \right.
\end{equation*}

If $j=m_{k-1}+1$ and $I\in \mathcal{F}_j$, then
$I$ contains exactly $2^{k-1}$ component intervals in $\mathcal{F}_{m_k}$ and $2^{k-1}-1$ gaps. Furthermore,
\begin{equation*}
|I|=\left\{
       \begin{array}{ll}
       (2^{k}-1)\delta_k & \text{ if $G_{l_k}\cap I=\emptyset$}, \\
       2^k\delta_k+k(\sqrt{2})^{k}\delta_k& \text{ if $G_{l_k}\subset I$}.
        \end{array}
     \right.
\end{equation*}

Inductively, suppose $m_{k-1}\leq j<m_{k}$ and $I\in \mathcal{F}_j$, $I$ contains exactly $2^{k-j+m_{k-1}}$ component intervals in $\mathcal{F}_{m_k}$ and $2^{k-j+m_{k-1}}-1$ gaps between them. Then
\begin{equation}\label{equation-0}
|I|=\left\{
       \begin{array}{ll}
       (2^{k-j+m_{k-1}+1}-1)\delta_k & \text{ if $G_{l_k}\cap I=\emptyset$},\\
       2^{k-j+m_{k-1}+1}\delta_k+k(\sqrt{2})^{k}\delta_k& \text{ if $G_{l_k}\subset I$}.
        \end{array}
     \right.
\end{equation}
Furthermore,
\begin{equation}\label{equation-1}
C(I,1)=\left\{
       \begin{array}{ll}
       |I|-\delta_k & \text{ if $j\neq m_{k-1}+[\frac{k}{2}]$},\\
       |I|-2\delta_k-k(\sqrt{2})^{k}\delta_k& \text{ if $j=m_{k-1}+[\frac{k}{2}], G_{l_k}\subset I $},\\
       |I|-\delta_k& \text{ if $j=m_{k-1}+[\frac{k}{2}], G_{l_k}\not\subset I$}.
        \end{array}
     \right.
\end{equation}

Now we consider the elements in $\mathcal{F}_{m_{k-1}+[\frac{k}{2}]}$. By the construction of $l_k$, there is a basic interval $I\in \mathcal{F}_{m_{k-1}+[\frac{k}{2}]}$ such that $G_{l_k}\subset I$, and $G_{l_k}$ is a gap in the $\mathcal{F}_{m_{k-1}+[\frac{k}{2}]+1}$. Furthermore, for $I\in \mathcal{F}_{m_{k-1}+[\frac{k}{2}]}$,
\begin{equation*}|I|=\left\{
\begin{array}{ll}
2^{k-[\frac{k}{2}]+1}\delta_k+k(\sqrt{2})^{k}\delta_k & $ if $ G_{l_k}\subset I,\\
(2^{k-[\frac{k}{2}]+1}-1)\delta_k & $ if $G_{l_k}\cap I=\emptyset.
\end{array}
\right.
\end{equation*}

Since
$$
1\leq\frac{2^{k-[\frac{k}{2}]+1}\delta_k+k(\sqrt{2})^{k}\delta_k}{(2^{k-[\frac{k}{2}]+1}-1)\delta_k}\rightarrow\infty
$$
as $k\rightarrow\infty$,
we find that the condition $(c)$ in Theorem \ref{theorem-3} does not hold for this set $E$.

Secondly, we prove this set is minimal by checking whether it satisfies the conditions in Theorem \ref{theorem-2} or not.

Condition $(\tilde{a})$: For sufficient large $j$, suppose $m_{j-1}\leq k<m_{j}$. Obviously, if $k\neq m_{j-1}+[\frac{j}{2}]$, then $I$ contains exactly a gap with length $\delta_j$ for any $I\in \mathcal{F}_k$, thus we have
\begin{equation}
\begin{split}
\beta_k&=\max\{\frac{\delta_j}{2^{j+m_{j-1}-k+1}\delta_j-\delta_j},
\frac{\delta_j}{2^{j+m_{j-1}-k+1}\delta_j+j(\sqrt{2})^{j}\delta_j}\}\\
&\leq\frac{\delta_j}{2^{j+m_{j-1}-k}\delta_j}=\frac{1}{2^{m_j-k}}.\nonumber
\end{split}
\end{equation}
If $k=m_{j-1}+[\frac{j}{2}]$, we have $\beta_k<1$ by the definition of $\beta_k$.

In conclusion, we have for sufficient large $j$,
\begin{equation}
\  \sum_{k=m_{j-1}}^{m_j-1}\beta_k=\sum_{k=\frac{j(j-1)}{2}}^{\frac{j(j+1)}{2}-1}\beta_k\leq2,\nonumber
\end{equation}
which implies that
$$
\lim\limits_{m\rightarrow\infty}\frac{1}{m}\sum_{j=0}^{m-1}\beta_{j}=0.
$$
That is $E$ satisfies Condition $(\tilde{a})$ of Theorem \ref{theorem-2}.

Condition $(\tilde{b})$: For sufficient large $k$, let $m_{k-1}\leq j<m_{k}$. If $j=m_{k-1}+[\frac{k}{2}]$ and $I\in \mathcal{F}_{j}$,
then $C(I,1)=2(2^{k-[\frac{k}{2}]}-1)\delta_k$ and
\begin{equation}
\begin{split}
 \Gamma_j&=\frac{2(2^{k-[\frac{k}{2}]}-1)\delta_k}{2(2^{k-[\frac{k}{2}]}-1)\delta_k+2\delta_k+k(\sqrt{2})^{k}\delta_k}\\
         &\geq\frac{2-\sqrt{2}}{4+k}.\nonumber
\end{split}
\end{equation}
If $j\neq m_{k-1}+[\frac{k}{2}]$ and $I\in \mathcal{F}_{j}$, from equations (\ref{equation-0}) and (\ref{equation-1}), we have $C(I,1)=|I|-\delta_k$ and
\begin{equation}
\Gamma_j=1-\frac{\delta_k}{|I|}\geq 1-\frac{1}{2^{k+m_{k-1}-j}}=1-\frac{1}{2^{m_{k}-j}}.\nonumber
\end{equation}
That is we have for sufficient large $k$,
\begin{equation}
\  \sum_{j=m_{k-1}}^{m_k-1}\log\Gamma_j\geq \sum_{j=m_{k-1}}^{m_k-1}\log(1-\frac{1}{2^{m_k-j}})-\log(4+k)+\log(2-\sqrt{2})\geq C-\log(4+k)\nonumber,
\end{equation}
where $C=\Sigma_{k=1}^{\infty}\log(1-\frac{1}{2^{k}})+\log(2-\sqrt{2})$.
For every $m>1$, there is a $m_0$ such that $m_{m_0-1}\leq m-1<m_{m_0}$. Thus
$$
\frac{1}{m}\sum_{j=0}^{m-1}\log\Gamma_j\geq \frac{1}{m}\sum_{k=1}^{m_0}\sum_{j=m_{k-1}}^{m_k-1}\log\Gamma_j\geq \frac{1}{m}(m_0C-\sum_{k=1}^{m_0}\log(4+k)).
$$
Note that $m_{m_0}=1+2+\cdots+m_0=\frac{m_0(m_0+1)}{2}$, we have
\begin{equation}
 \lim\limits_{m\rightarrow\infty}\frac{\sum_{k=1}^{m_0}\log(4+k)}{m}=0.\nonumber
\end{equation}
In conclusion, we have
$$
0\geq\lim_{m\rightarrow\infty}\frac{1}{m}\sum_{j=0}^{m-1}\log\Gamma_j\geq
\lim_{m\rightarrow\infty}\frac{1}{m}(m_0C-\sum_{k=1}^{m_0}\log(4+k))=0.
$$

Condition $(\tilde{c})$: For sufficient large $k$, let $m_{k-1}< j\leq m_{k}$. If $m_{k-1}+[\frac{k}{2}]+2\leq j\leq m_k$, for any $I\in\mathcal{F}_{j}$, they have the same length, and $|Fa(I)|$ have the same length too. That is for $m_{k-1}+[\frac{k}{2}]+2\leq j\leq m_k$, we have
\begin{equation}
\  \gamma_j=\frac{|I|}{|Fa(I)|}=\frac{(2^{k+m_{k-1}-j}+2^{k+m_{k-1}-j}-1)\delta_k}{(2^{k+m_{k-1}-j+1}+2^{k+m_{k-1}-j+1}-1)\delta_k}<\frac{1}{2}.\nonumber
\end{equation}
Set $\alpha=\frac{1}{2}$, it is easy to see that condition $(\tilde{c})$ holds for $E$.

\begin{remark}With some changes of this example, for example, we let $n_k=2^{k}$ and $c_k\cdot 2^{k+1}=\frac{2}{3}$. For $k\geq1$, let $\eta_{k,0}=\eta_{k,n_k}=0$. For $k\geq2$, let $\eta_{k,1}=\eta_{k,2}=\cdots=\eta_{k,n_k-2}=\delta_k=c_1c_2\cdots c_k$ and $\eta_{k,n_k-1}=\frac{1}{3}\delta_{k-1}+2\delta_k$. We have constructed a homogeneous perfect set which satisfies the condition $(3)$ in Theorem \ref{Hdim} and has Hausdorff dimension 1. But it does not satisfy some condition in Theorem \ref{theorem-2}. So our main theorem does not supply a general method to deal with homogeneous perfect set which satisfies the condition $(3)$ of Theorem \ref{Hdim}.
\end{remark}

\noindent\textbf{Acknowledgements.} This work was supported by National Natural Science Foundation of China (Grant Nos.\,11301165, 11371126, 11571099). The authors wish to thank Prof. Jiang Yueping for his support and encouragement of this work.

\end{document}